\documentclass[12pt,letterpaper]{article}
\usepackage{amsfonts,amssymb,latexsym}

\newtheorem{theorem}{Theorem}[section]
\newtheorem{lemma}{Lemma}[section]

\newtheorem{definition}{Definition}[section]
\newtheorem{corollary}{Corollary}[section]

\newtheorem{remark}{Remark}[section]

\newenvironment{itemizz}{\begin{itemize}\setlength{\itemsep}{-1mm}} %
{\end{itemize}}

\newcommand{\eop}{$\quad\spadesuit$}
\newenvironment{proof}{{\bf Proof.}}{\eop\medskip}

\begin{document}

\author{Ramiro de la Vega}
\title{Countable Tightness, Elementary Submodels and Homogeneity}
\maketitle

\begin{abstract}
We show (in $ZFC$) that the cardinality of a compact homogeneous space of countable tightness is no more than $2^{\aleph_0}$.
\end{abstract}

\section{Introduction}

A space $X$ is said to have \textit{countable tightness} ($t(X)=\aleph_0$) if whenever $A\subseteq X$ and $x\in \overline{A}$, there is a countable $B \subseteq A$ such that $x \in \overline{A}$. A space $X$ is \textit{homogeneous} if for every $x,y\in X$ there is a homeomorphism $f$ of $X$ onto $X$ with $f(x) =y$. It is known (see \cite{juh}) that any compact space of countable tightness contains a point with character at most $2^{\aleph_0}$; if the space is also homogeneous then it follows that $|X|\leq 2^{2^{\aleph_0}}$. In \cite{arc}, Arhangel'skii conjectured that in fact $|X|\leq 2^{\aleph_0}$ for any such space. The main goal of this paper is to give a proof (in $ZFC$) to Arhangel'skii's conjecture. This is achieved in Theorem \ref{goal}. As a corollary of this we also confirm a conjecture of I. Juh\'asz, P. Nyikos and Z. Szentmikl\'ossy (see \cite{nyi}), stating that it is consistent that every homogeneous $T_5$ compactum is first countable.

Our main tool will be the ``Elementary Submodels technique": Given a topological space $(X,\tau)$, let $M$ be an elementary submodel of $H(\theta)$ (i.e. the set of all sets of hereditary cardinality less than $\theta$) for a ``large enough" regular cardinal $\theta$. Then one uses properties of $X \cap M$, $\tau \cap M$ and $M$ itself to get results about $(X,\tau)$. A model $M$ is \textit{countably closed} if any countable sequence of elements of $M$ is in $M$ (i.e. $M^\omega \subseteq M$). For more details and a good introduction to the technique see \cite{dow}. Let us just say that in each specific application, one takes $\theta$ large enough for $H(\theta)$ to contain all sets of interest in the context under discussion. In this sense we will just say $M \prec {\bf V}$. In Section \ref{models} we prove some basic facts in the context of elementary submodels of countably tight compact spaces; we also give answer (Theorem \ref{normal}) to a question of L.R. Junqueira and F. Tall.

\section{Elementary submodels}\label{models}

Fix a compact Hausdorff space $\left(X,\tau \right)$ with $t(X)=\aleph_0$ and fix a countably closed $M\prec {\bf V}$ with $X,\tau \in M$. Let $Z=\overline{X\cap M}\subseteq X$ with the subspace topology.

One of the main goals of this section is to show that $Z$ is a retract of $X$. The following result suggests what the retraction will be.

\begin{lemma}\label{exists}
For every $x \in X$ there is a $q_x \in Z$ such that for all $U \in \tau \cap M$ either $q_x \notin U$ or $x \in U$.
\end{lemma}

\begin{proof} 
Fix $x \in X$ and assume there is not such a $q_x$. Then for each $q \in Z$ we can fix a $U_q \in \tau \cap M$ such that $q \in U_q$ and $x \notin U_q$. Since $Z$ is compact we get that $Z\subseteq \bigcup_{q \in Q}U_q$ for some finite $Q\subseteq Z$. On the other hand $x \notin \bigcup_{q \in Q}U_q \in M$ so by elementarity there is an $x' \in (X \cap M) \setminus \bigcup_{q \in Q}U_q$ which is impossible.
\end{proof}

Just by elementarity and the fact that $X$ is Hausdorff it is immediate that $\tau \cap M$ separates points in $X \cap M$. We prove now that in fact $\tau \cap M$ separates points in $Z$.

\begin{lemma}\label{separates}
If $p_0,p_1 \in Z$ and $p_0 \ne p_1$ then there are $U_0,U_1 \in \tau \cap M$ such that $p_0 \in U_0$, $p_1 \in U_1$ and $U_0 \cap U_1=\emptyset$.
\end{lemma}

\begin{proof}
Given $p_0,p_1 \in Z$ with $p_0 \ne p_1$ we use regularity of $X$ to fix $U'_i,V_i \in \tau$ for $i \in 2$ such that $p_i \in V_i$, $\overline{V}_i\subseteq U'_i$ and $U'_0 \cap U'_1 = \emptyset$. Since $t(X)=\aleph_0$, there are countable $A_0,A_1 \subseteq M$ such that $A_i \subseteq V_i$ and $p_i \in \overline{A}_i$ for $i \in 2$. Since $M$ is countably closed we have that $A_0,A_1 \in M$ and hence by elementarity $M \models \exists U_0,U_1 \in \tau \left[\ \overline{A}_0 \subseteq U_0, \overline{A}_1 \subseteq U_1, U_0 \cap U_1 = \emptyset \ \right]$.
\end{proof}

\begin{corollary}\label{unique}
For every $x \in X$ there is a unique $q_x \in Z$ such that for all $U \in \tau \cap M$ either $q_x \notin U$ or $x \in U$.
\end{corollary}

\begin{proof}
Immediate from Lemmas \ref{exists} and \ref{separates}.
\end{proof}

\begin{definition} In view of the last result we define the function $r_M:X \rightarrow Z$ by $r_M(x)=q_x$ for $x \in X$.
\end{definition} 

\begin{lemma}\label{continuous}
$r_M$ is continuous.
\end{lemma}

\begin{proof}
Fix $W \in \tau$ (not necessarily in $M$) with $W \cap Z \ne \emptyset$, fix $x \in r_M^{-1}(W)$ and let $q=r_M(x)$. We need to show that there is a $V \in \tau$ such that $x \in V \subseteq r_M^{-1}(W)$.

For each $p \in Z \setminus W$ we use Lemma \ref{separates} to get $U_p,V_p \in \tau \cap M$ such that $p \in U_p$, $q \in V_p$ and $U_p \cap V_p = \emptyset$. Since $Z \setminus W$ is closed (and hence compact) we get that $Z \setminus W \subseteq U:=\bigcup_{p \in P}U_p$ for some finite $P\subseteq Z \setminus W$. Clearly $q \in V:=\bigcap_{p \in P}V_p$ and $U \cap V = \emptyset$. Also by elementarity we have that $U,V \in M$.

Since $q \in V \in M$ we get (from the definition of $r_M$) that $x \in V$. To see that $V \subseteq r_M^{-1}(W)$, fix $y \in V$ and note that if $r_M(y)$ was in $Z \setminus W$ then it would be in $U$ and since $U \in M$ we would get that $y \in U$; but this is impossible since $U \cap V = \emptyset$ and hence $r_M(y) \in W$.
\end{proof}

It is clear that $r_M(x)=x$ for all $x \in Z$. So we actually have that $r_M$ is a retraction. Looking closer at the last proof, we see that $q \in V \cap Z \subseteq W \cap Z$, so we also showed the following

\begin{corollary}\label{base}
The set $\{U \cap Z : U \in \tau \cap M\}$ is a base for the topology of $Z$.
\end{corollary}

For reference, we summarize our results in the following

\begin{theorem}\label{main}
Let $(X,\tau)$ be a compact space of countable tightness and let $M\prec {\bf V}$ be countably closed with $X,\tau \in M$. Then $r_M:X\rightarrow \overline{X\cap M}$ is a retraction and $\{U \cap \overline{X\cap M} : U \in \tau \cap M\}$ is a base for $\overline{X\cap M}$.
\end{theorem}

In \cite{tall}, Junqueira and Tall define the space $X_M$ as the set $X \cap M$ with the topology $\tau_M$ generated by $\{U \cap M : U \in \tau \cap M\}$. They ask (Problem 7.22) if there is a consistent example of a compat $T_2$ space $X$ with countable tightness and a countably closed $M$ such that $X_M$ is not normal. A consequence of Theorem \ref{main} is that $X_M$ is a subspace of $X$; this was already proved in \cite{tall} (Theorem 2.11). In the next theorem we make use of this fact to give a negative answer to their question.

\begin{theorem}\label{normal}
If $(X,\tau)$ is a compact Hausdorff space of countable tightness and $M$ is countably closed then $X \cap M$ ($=X_M$) is normal.
\end{theorem}

\begin{proof}
Fix two disjoint closed $E,F \subseteq X \cap M$. We claim that $\overline{E}$ and $\overline{F}$ are still disjoint. Therefore, since $\overline{X \cap M}$ is compact (and hence normal), $\overline{E}$ and $\overline{F}$ (and hence $E$ and $F$) can be separated.

Now suppose (seeking a contradiction) that $p \in \overline{E} \cap \overline{F}$. Since $X$ is countably tight, there are countable $C_E \subseteq E$ and $C_F \subseteq F$ such that $p \in \overline{C_E} \cap \overline{C_F}$. Since $C_E$ and $C_F$ are countable, they are in $M$ and thus by elementarity there is a $q \in X \cap M$ with $q \in \overline{C_E} \cap \overline{C_F}$. But this is impossible since $E$ and $F$ are closed and disjoint in $X \cap M$. Hence $\overline{E} \cap \overline{F} = \emptyset$.
\end{proof}

\section{Homogeneity}

The following lemma was proved in \cite{juh}.

\begin{lemma} If $X$ is a compact $T_2$ space of countable tightness then there are a countable set $S\subseteq X$ and a non-empty closed $G_\delta$ set $H$ in $X$ with $H \subseteq \overline{S}$.
\end{lemma}

In general, given a point $p \in X$ one cannot expect to get $p \in H$ in the previous lemma. For example if $X=\kappa \cup \{\infty\}$ is the one point compactification of an uncountable discrete $\kappa$ then $p=\infty$ is a counter example. In this space the clousure of any countable set is countable, but on the other hand any $G_\delta$ subset containing $p$ must be uncountable. However the following is obvious.

\begin{corollary}\label{delta} If $X$ is a compact homogeneous $T_2$ space of countable tightness and $p\in X$ then there are a countable set $S\subseteq X$ and a closed $G_\delta$ set $H$ in $X$ with $p\in H$ and $H \subseteq \overline{S}$.
\end{corollary}

Now we are ready to prove our first important result.

\begin{theorem}\label{tight}
Suppose $(X,\tau)$ is a compact homogeneous Hausdorff space of countable tightness. Then $w(X)\leq 2^{\aleph_0}$.
\end{theorem}

\begin{proof}
By Corollary \ref{delta} and regularity of $X$ we can fix functions $\psi:X\times \omega \rightarrow \tau$ and $S:X \rightarrow [X]^{\aleph_0}$ such that for all $x\in X$ and $n\in \omega$:
\begin{itemizz}
\item[1.]
$x\in \psi(x,n)$ and $\overline{\psi(x,n+1)}\subseteq \psi(x,n)$.
\item[2.]
$H(x):=\bigcap_{m\in \omega}\psi(x,m) = \bigcap_{m\in \omega}\overline{\psi(x,m)} \subseteq \overline{S(x)}$.
\end{itemizz} 

Now fix a countably closed $M\prec {\bf V}$ with $X,\tau,\psi,S \in M$ and $|M|\leq 2^{\aleph_0}$. Let $Z=\overline{X\cap M}$ and $r_M:X\rightarrow Z$ as in Section \ref{models}. We know (by Theorem \ref{main}) that $w(Z)\leq |\tau \cap M| \leq 2^{\aleph_0}$. We will show that in fact $X=Z$ and hence $w(X)\leq 2^{\aleph_0}$.

Fix now $x \in X$ and let $p = r_M(x) \in Z$. By Theorem \ref{main}, for each $n \in \omega$ we can find $U_n \in \tau \cap M$ such that $\overline{\psi(p,n+1)} \cap Z \subseteq U_n \cap Z \subseteq \psi(p,n) \cap Z$. Also by countable tightness of $X$ we can fix $A \in [X \cap M]^{\aleph_0}$ such that $p \in \overline{A}$. Now since $A$, $\psi$, $H$ and $<U_n:n \in \omega>$ are all in $M$ (since $M$ is countably closed), we get by elementarity that there must be a $q \in X \cap M$ such that (note that the three conditions are true of $p$ in ${\bf V}$):
\begin{itemizz}
\item[1.]
$q\in \bigcap_{n \in \omega}U_n$.
\item[2.]
$q \in \overline{A}$.
\item[3.]
$\overline{A} \cap \bigcap_{n \in \omega}U_n = \overline{A} \cap H(q)$.
\end{itemizz} 

Since $p$ satisfies conditions 1 and 2, we get that $p \in H(q)$ and hence $x \in H(q)$ by definition of $r_M$. But on the other hand $S(q)\subseteq X \cap M$ and $H(q) \subseteq \overline{S(q)}$ and therefore $x \in Z$. This shows that $X = Z$ which is what we wanted.
\end{proof}

\begin{remark}
With more work one can show that in fact $X \subseteq M$, showing then that $|X|\leq 2^{\aleph_0}$. However we shall use a different (perhaps simpler) strategy to get this result. 
\end{remark}

The following is a well known fact about compact homogeneous spaces. In fact it was shown in \cite{mill} that it also holds for any compact power homogeneous space.

\begin{lemma}
If $X$ is compact and homogeneous, then $|X|\leq w(X)^{\pi\chi(X)}$.
\end{lemma}

Putting together the two previous results and the fact that $\pi\chi(X) \leq t(X)$ for any compact space, we immediately get our main result.

\begin{theorem}\label{goal}
If $X$ is compact, homogeneous and $t(X)=\aleph_0$ then $|X|\leq 2^{\aleph_0}$.
\end{theorem}

Using the fact that $|X|=2^{\chi(X)}$ for compact homogeneous spaces, we get:

\begin{corollary}
$(2^{\aleph_0}<2^{\aleph_1})$ If $X$ is compact, homogeneous and $t(X)=\aleph_0$ then $X$ is first countable.
\end{corollary}

In \cite{mill}, J. van Mill asked whether every $T_5$ (i.e. hereditarily normal) homogeneous compact space has cardinality $2^{\aleph_0}$. In \cite{nyi}, I. Juh\'asz, P. Nyikos and Z. Szentmikl\'ossy proved that the answer is yes in forcing extension resulting by adding $(2^{\aleph_1})^{\bf V}$ Cohen reals. They also showed that after adding $\aleph_2$ Cohen reals, every $T_5$ homogeneous compact space has countable tightness. Putting this together with Theorem \ref{goal} and assuming for example that $2^{\aleph_0}=\aleph_2$ and $2^{\aleph_1}=\aleph_3$ in the ground model, we get a confirmation of a conjecture proposed by them in the same paper.

\begin{theorem}
It is consistent that every homogeneous $T_5$ compact space is first countable.
\end{theorem}

\end{document}